\documentclass[11pt]{amsart}
\usepackage{amscd}

\newcommand{\ncom}{\newcommand}
\ncom{\ra}{\rightarrow}
\ncom{\lra}{\longrightarrow}
\ncom{\Rar}{\Rightarrow}
\ncom{\Lra}{\Leftrightarrow}
\ncom{\noin}{\noindent}

\ncom{\socle}{\operatorname{socle}}
\ncom{\type}{\operatorname{type}}
\ncom{\depth}{\operatorname{depth}}
\newtheorem{thm}{Theorem}[section]
\newtheorem{lemma}[thm]{Lemma}
\newtheorem{cor}[thm]{Corollary}
\newtheorem{pro}[thm]{Proposition}
\newtheorem{example}[thm]{Example}
\newtheorem{definition}[thm]{Definition}
\newtheorem{remark}[thm]{Remark}
\newtheorem{notation}[thm]{Notation}
\newtheorem{blank}[thm]{}

\ncom{\sz}{\scriptsize}
\ncom{\bib}{\bibitem}
\ncom{\ol}{\overline}

\newcommand{\ZZ}{\mathbb Z}
\newcommand{\NN}{\mathbb N}
\newcommand{\ga}{\alpha}

\def\m{{\mathfrak m}}
\def\n{{\mathfrak n}}

\begin{document}

\title{ On  fiber cones of $\m$-primary ideals }

\author{  A. V. Jayanthan}

\address{Harish-Chandra Research Institute, Chhatnag Road, Jhusi,
Allahabad 211019, INDIA.}
\email{jayan@mri.ernet.in}

\author{ Tony J. Puthenpurakal}
\address{ Department of Mathematics, Indian Institute of Technology
Bombay, Powai, Mumbai 40076, INDIA.}
\email{tputhen@math.iitb.ac.in}

\author{ J. K. Verma}
\address{ Department of Mathematics, Indian Institute of Technology
Bombay, Powai, Mumbai 40076, INDIA.}
\email{jkv@math.iitb.ac.in}
\thanks{AMS Classification 2000: 13H10, 13H15, 13A30 (Primary) 13C15,
13A02 (Secondary)}
\keywords{ fiber cones, mixed multiplicities, joint reductions,
Cohen-Macaulay fiber cones, Gorenstein fiber cones, 
ideals having minimal and almost minimal mixed multiplicities.}

\begin{abstract} Two formulas for the multiplicity of the fiber cone
  $F(I)=\oplus_{n=0}^{\infty} I^n/\m I^n$ of an $\m$-primary ideal of
  a $d$-dimensional Cohen-Macaulay local ring $(R,\m)$ are derived  in
  terms of the mixed multiplicity $e_{d-1}(\m | I),$ the multiplicity
  $e(I)$ and superficial elements. As a consequence, the
  Cohen-Macaulay property of $F(I)$ when $I$ has minimal mixed
  multiplicity or almost minimal mixed multiplicity is characterized
  in terms of reduction number of $I$ and lengths of certain ideals.
  We  also characterize Cohen-Macaulay and Gorenstein property of
  fiber cones of  $\m$-primary ideals  with a $d$-generated minimal
  reduction $J$ satisfying  (i) $\ell(I^2/JI)=1$ or (ii)
  $\ell(I\m/J\m)=1.$  
\end{abstract}

\maketitle

\thispagestyle{empty}

\section{\bf Introduction}
The objective of this paper is to study Cohen-Macaulay and Gorenstein
properties of the fiber cone $F(I)=\oplus_{n=0}^{\infty} I^n/\m I^n$
of an $\m$-primary ideal $I$ of a $d$-dimensional Cohen-Macaulay local
ring $(R,\m)$ in terms of the invariants such as the multiplicity
$e(I),$ the mixed multiplicity  $e_{d-1}(\m | I)$ and
reduction number of $I.$  

In order to state the main results, we recall  necessary
definitions first.  Let $I$ be an $\m$-primary ideal of a
$d$-dimensional local ring $(R, \m).$ 
The Hilbert function  $HF(F(I),n)$ of
the fiber cone $F(I)$ is defined as $HF(F(I),n)= \ell (I^n/\m I^n)$, 
where $\ell$ denotes the length function.  The function $HF(F(I),n)$ is a
polynomial $HP(F(I),n)$ in $n$ of degree $d-1$ for all large $n.$ We
write this polynomial as
$$ HP(F(I),n)
=f_0(I)\binom{n+d-1}{d-1}-f_1(I)\binom{n+d-2}{d-2} 
 + \cdots +(-1)^{d-1} f_{d-1}(I),$$
for certain integers $f_0(I), f_1(I), \ldots, f_{d-1}(I).$ The number
$f_0(I)$ is called the {\em multiplicity} of $F(I).$

\noindent
{\bf Multiplicities  and reductions:} 
For an $\m$-primary ideal $I$ in a Noetherian
local ring $R$ of dimension $d$, let $HF(I,n):= \ell(R/I^n)$ denote
the {\it Hilbert-Samuel function} of $I$. It is well known that this
function coincides with a polynomial $HP(I,n)$ of degree $d$. Write
the polynomial as:
$$
HP(I,n) = e_0(I) {n+d-1 \choose d} - e_1(I){n+d-2 \choose d-1} +
\cdots + (-1)^de_d(I).
$$
The coefficient $e_0(I),$ also denoted as $e(I),$  is called the 
multiplicity of $I$. 

Now we recall some basic facts about {\em reductions} from
\cite{nr}.  An ideal $K \subseteq I$ is called a {\em reduction } of $I$
if there exists a nonnegative integer $n$ such that $KI^n=I^{n+1}.$ If
$K$ is minimal with respect to inclusion among reductions of $I,$ then
it is called a {\em minimal reduction} of $I.$ The {\em reduction
number} $r(I)$ of $I$ is the least integer $n$ such that
$JI^n=I^{n+1},$ where $J$ varies over all minimal reductions of $I.$
If $R/\m$ is infinite, then all minimal reductions of $I$ are
generated by the same number of elements called the {\em analytic spread}
of $I.$ The analytic spread  of $I,$ is the Krull dimension of
the fiber cone $F(I).$ It is easy to see that if $J$ is a reduction of $I,$ 
then $e(I)=e(J).$ 

\medskip
\noindent
{\bf Mixed multiplicities and joint reductions:} 
Mixed multiplicities and joint reductions of ideals are analogues of 
reductions and multiplicities of ideals.  
Let $I_1, I_2, \ldots, I_r$  be $\m$-primary ideals.  The {\em Bhattacharya
function} of  $I_1, I_2, \ldots, I_r$   is the numerical function 
$BF(n_1, n_2, \ldots, n_r):\NN^r \lra
\NN,$ defined by $BF(n_1, n_2, \ldots, n_r)  
=\ell(R/I_1^{n_1} I_2^{n_2}\ldots I_r^{n_r}).$ By \cite{te}, for all  
$n_1, n_2, \ldots, n_r ,$ large, the  Bhattacharya function 
 is given by a polynomial $BP(n_1, n_2, \ldots, n_r)$ of
total degree $d$ in  $n_1, n_2, \ldots, n_r.$ For 
$\ga = ( \ga_1, \ga_2, \ldots, \ga_r) \in \NN^r$ we put  
$| \ga |=\ga_1+\ga_2 + \cdots + \ga_r.$
We write the  Bhattacharya polynomial in the form 
$$
BP(n_1, n_2, \ldots, n_r)=\sum_{ |\ga | \leq d  }
e_{\ga}
\binom{n_1+\ga_1}{\ga_1}  \binom{n_2+\ga_2}{\ga_2}  \cdots 
\binom{n_r+\ga_r}{\ga_r},
$$
where $e_{\ga}$ are certain integers. 
Let $\ga=(\ga_1, \ga_2, \ldots, \ga_r) \in \NN^r$ and $|\ga|=d.$
A  multiset consisting of 
$\ga_1$ copies of $I_1, \ga_2$ copies of 
$I_2,  \ldots , \ga_r$ copies of $I_r,$ will be denoted by 
$(I_1^{[\ga_1]}| I_2^{[\ga_2]}| \cdots | I_r^{[\ga_r]}).$ In case 
$|\ga|=d,$ we write   
$e_{\ga}=e_{\ga}(I_1^{[\ga_1]}| I_2^{[\ga_2]}| 
\cdots | I_r^{[\ga_r]}).$
These  integers are positive and  are called {\em mixed multiplicities} 
of the ideals $I_1, I_2, \ldots, I_r.$ 
When $r=2,$ and $i+j=d,$ we adopt the simpler notatioon 
$e_{(i,j)}(I^{[i]}|J^{[j]})=e_j(I|J).$  
D. Rees proved in \cite{r1} that $e_0(I|J)=e(I)$ and $e_d(I|J)=e(J).$

D. Rees  introduced   joint reductions
in  \cite{re} for calculating mixed multiplicities.  
Let $(R,\m)$ be a $d$-dimensional local ring. Let $I_1, I_2, \ldots, I_d$ 
be $\m$-primary ideals. We say that 
 $a_1 \in I_1, a_2 \in I_2, \ldots, a_d \in I_d $ is a  
 \emph{joint reduction} of $I_1, I_2, \ldots, I_d$ if 
$a_1 I_2I_3\ldots I_d+a_2I_1I_3\ldots I_d+ \cdots+a_dI_1I_2 \ldots I_{d-1}$
is a reduction of $I_1I_2\ldots I_d.$ 
D. Rees showed in Theorem 2.4 of \cite{re} that 
if $(a_1, a_2, \ldots, a_d)$ is a joint reduction of 
the multiset $(I_1^{[\ga_1]}| I_2^{[\ga_2]}| \cdots | I_r^{[\ga_r]}),$ 
where $|\ga|=d,$ then
$$e((a_1, a_2, \ldots, a_d))
=e_{\ga}(I_1^{[\ga_1]} |  I_2^{[\ga_2]} |  \cdots | I_r^{[\ga_r]}).$$

\medskip
\noindent
{\bf Rings and ideals of minimal and almost minimal multiplicity:} Let $\mu(I)$ denote the minimum number of elements required to
generate an ideal $I.$ For a Cohen-Macaulay local
ring $(R,\m)$ of dimension $d$, $e(\m) \geq \mu(\m)-d+1$. A Cohen-Macaulay 
local  ring is
said to have minimal multiplicity (resp. almost minimal multiplicity)
if $e(\m) = \mu(\m)-d+1$ (resp. $e(\m) = \mu(\m)-d+2$). J. D. Sally
studied  Cohen-Macaulay local rings  of minimal and almost minimal 
multiplicity. She  proved 
that the associated graded ring $G(\m) := \oplus_{n\geq 0}\m^n/\m^{n+1}$ is
Cohen-Macaulay  when $R$ is Cohen-Macaulay with minimal multiplicity 
\cite{s1}. She conjectured that if the
ring has almost minimal multiplicity, then $G(\m)$ has depth at least
$d-1$, \cite{s2}. This conjecture was proved independently by Wang
\cite{w} and M. E. Rossi and G. Valla \cite{rv}. Later on M. E. Rossi
generalized the conjecture of J. D. Sally to the case of $\m$-primary
ideals. She proved that if $e(I) = \ell(I/I^2) + (1-d)\ell(R/I) + 1$,
then $\depth G(I) \geq d-1$ \cite{r}. It is easy to see that $e(I) =
\ell(I/I^2) + (1-d)\ell(R/I) + 1$ if and only if for any minimal
reduction $J$ of $I,$ $\ell(I^2/JI)=1.$  

\begin{definition}
An $\m$-primary ideal $I$ of a Cohen-Macaulay local ring  satisfying
the condition $\ell(I^2/JI)=1,$ for any minimal reduction J, is
called a Sally ideal.
\end{definition}

The notions of minimal multiplicity and almost minimal multiplicity have
been generalized in many directions. It was proved  in  \cite{drv}
that for an $\m$-primary ideal $I$ of a Cohen-Macaulay local ring $(R,
\m),$ $e_{d-1}(\m|I) \geq \mu(I)-d+1.$ We say that $I$ has {\em
minimal mixed multiplicity} if $e_{d-1}(\m|I)= \mu(I)-d+1$ and $I$ has
{\em almost minimal mixed multiplicity} if $e_{d-1}(\m | I)=
\mu(I)-d+2.$ In \cite{drv} and \cite{dv}, the Cohen-Macaulay property of
fiber cones of ideals with minimal and almost minimal mixed
multiplicities was studied.  J. Chuai, in \cite{c}, proved that for
an $\m$-primary ideal $I$ in a Cohen-Macaulay local ring $(R,\m)$,
$e(I) \geq \mu(I) - d + \ell(R/I)$. S. Goto termed an ideal to have
minimal multiplicity if $e(I) = \mu(I) -d+\ell(R/I)$, \cite{g2}. He
studied many properties of the associated graded ring, the fiber cone
and the Rees algebra of ideals with minimal multiplicity. In
\cite{jv2}, fiber cones of ideals having almost minimal multiplicity
are studied,  i.e. ideals with the property $e(I) = \mu(I) -
d+\ell(R/I) + 1$.

\medskip
\noindent
{\bf Main results:} In this paper, we consider the Cohen-Macaulay and 
Gorenstein properties of fiber cones of Sally ideals, ideals with minimal and
almost minimal (mixed) multiplicity. We assume, in the rest of this
section, that {\em $(R, \m)$ is a $d$-dimensional Cohen-Macaulay local ring
with infinite residue field.} 

In section 2, we  obtain two  formulas for $f_0(I)$   in terms of
$e_{d-1}(\m | I), e(I)$ and superficial elements for $\m$ and $I$
in the sense of Rees. 

In Section 3, as a consequence of these  formulas we recover  one of
the main results of \cite {drv} to the effect that for an ideal $I$ of
minimal mixed multiplicity, $F(I)$ is Cohen-Macaulay if and only if
$r(I) \leq 1.$ For an ideal $I$ of almost minimal mixed multiplicity,
we show that either $f_0(I)=e_{d-1} (\m | I)-1$ or
$f_0(I)=e_{d-1}(\m | I).$ In the former case, $F(I)$ is
Cohen-Macaulay if and only if $r(I)\leq 1$ and in the latter case,
$F(I)$ is Cohen-Macaulay if and only if $r(I)=2$ and $\ell(I^2/JI+\m
I^2)=1.$  In \cite{dv} this  result was proved under depth assumptions
on $G(I)$.  We improve  it by carefully using the  multiplicity
formula for the fiber cone. If $I$ is a Sally ideal with a minimal
reduction $J$, then we show that $F(I)$ is Cohen-Macaulay if and only
if $\m I^2=\m JI$ if and only if the Hilbert series of $F(I)$ is
$$HS(F(I),t) = \frac{1+(\mu(I)-d)t+t^2+\cdots+t^r}{(1-t)^d}.$$
In Section 4, we study the Gorenstein property of the
Cohen-Macaulay fiber cones. For this purpose we use Macaulay's theorem
about symmetry of the $h$-vector in the Hilbert series of $F(I)$ and
certain bilinear forms. It is fairly easy to show that for ideals of
reduction number $1,$ $F(I)$ is Gorenstein if and only if
$\mu(I)=d+1.$ We show that if  $r(I)=2,$ then $F(I)$ is
Gorenstein if and only if $$(I^2\m+JI:I) \cap I=\m I +J \quad
\text{and} \quad \ell(I^2/JI+\m I^2)=1.$$ For Sally ideals of
reduction number at least $3,$ we  show that $F(I)$ is Gorenstein
if and only if $\mu(I)=d+1.$

In Section 5, we characterize the Gorenstein property of fiber cones
of ideals of almost minimal multiplicity when $G(I)$ is
Cohen-Macaulay. We show that for such ideals $F(I)$ is Gorenstein
if and only if $I \cap (\m J:I)=\m I +J$ and when this is the case,
$\mu(I) \leq \mu(\m) +d.$  

In  Section 6, we illustrate our results with a few examples.

\vskip 2mm
\noindent
{\it Acknowledgments: } Part of this work was done while A. V. Jayanthan
 was visiting University of Genova. He is  thankful to the
University for the hospitality and A. Conca, M. E. Rossi and G. Valla
for  useful conversations. He also  thanks the  Indian Institute 
of Technology Bombay, Mumbai for its kind hospitality during his visit, 
while part of this work was done. J. K. Verma thanks the Abdus Salam 
International  Centre for Theoretical Physics, Trieste (Italy) for inviting 
him to visit the Centre as an {\em Associate Member} during May-June 2005, when 
this work was completed. The authors thank the referee for pertinent
comments.

\section{\bf Multiplicity formulas  for fiber cones}

Throughout this section $(R,\m)$ will denote a local ring.
In this section we derive two formulas for the multiplicity of the
fiber cone of an $\m$-primary ideal $I$ in $R$  in terms
of the mixed multiplicity $e_{d-1}(m|I)$ and the multiplicity $e(I).$
These formulas are in terms of superficial sequences for a set of ideals
in the sense of Rees.  We begin with a discussion of  
 superficial sequences and their relevance to joint reductions and hence 
mixed multiplicities.  As we need these for $\m$-primary ideals, 
we  restrict  our discussion to only such ideals. We begin by recalling 
the following definitions and results from I. Swanson's thesis \cite{sw}
for convenience of the reader.    

\begin{definition} 
[cf. Definition 1.14 of  \cite{sw}]  

Let $(R,\m)$ be a local ring. 
Let $I_1, I_2, \ldots, I_r$ 
be $\m$-primary $R$-ideals. An element $a \in I_1$ is called superficial 
for the ideals $I_1, I_2, \ldots, I_r $ if $\dim(R/(a))=\dim(R)-1$
and for some nonnegative integer $c$ and for all 
$n_1 > c, n_2,  \ldots,  n_r \geq 0,$
$$  \left(I_1^{n_1}I_2^{n_2}\ldots I_r^{n_r} : a \right)
  \cap I_1^cI_2^{n_2}\ldots I_r^{n_r} 
=  I_1^{n_1-1}I_2^{n_2}\ldots I_r^{n_r}.
$$ 
\end{definition} 

\begin{definition} 
[cf. Page 20 of  \cite{sw}]  
A sequence $a_1, a_2, \ldots, a_r$ of elements in $R$ is called 
a superficial sequence for the  ideals $I_1, I_2, \ldots, I_r$
if $a_i \in I_i$  and the image of $a_i$ in 
$R_{i-1}=R/(a_1, a_2, \ldots, a_{i-1})$ is superficial for the 
images of the ideals
$I_i, I_{i+1}, \ldots, I_r$ in $R_{i-1}$ for $i=1, 2, \ldots, r.$ 
    
\end{definition} 

\begin{thm} 
[cf. Theorem 1.16 of  \cite{sw}]   
\label{jr}
Let  $(R,\m)$ be  of positive dimension $d$ with   $R/\m$  infinite. 
Then superficial elements  exist. Moroever 
if ${\bf a}=a_1, a_2, \ldots, a_d$ is a superficial sequence for the 
$\m$-primary  ideals ${\bf I}=I_1, I_2, \ldots, I_d,$  then ${\bf a}$ 
is a joint reduction of ${\bf I}.$   
\end{thm} 

\noindent
Inspired by Rees' construction of joint reductions in his fundamental paper 
\cite {re} on joint reductions and mixed multiplicities, we introduce the 
following: 
\begin{definition} 
An element $a\in I_1$ is called  Rees-superficial  for the
$\m$-primary ideals     $I_1, I_2, \ldots, I_r$ if for all 
large $n_1$  and all nonnegative integers $n_2, n_3, \ldots, n_r,$ 
$$(a) \cap I_1^{n_1}  I_2^{n_2}  \ldots I_r^{n_r} =   
(a)I_1^{n_1-1}  I_2^{n_2}  \ldots I_r^{n_r}.$$
\end{definition}

\begin{definition}
A sequence $a_1, a_2, \ldots, a_r$ is called Rees-superficial 
for the ideals $I_1, I_2, \ldots, I_d,$ if the image of $a_{i}$ 
in $R_{i-1}=R/(a_1, a_2, \ldots, a_{i-1})$ is Rees-superficial for the 
images of  $I_{i}, \ldots, I_r$ in $R_{i-1} $ for $i=1, \ldots, r.$  
\end{definition}  

\begin{lemma}[{\bf Rees' Basic Lemma,} cf. Lemma 1.2 of \cite{re}] Let 
$I_1, I_2, \ldots, I_r$ 
be ideals of $R$ where  $R/\m$ is infinite.
Let ${\mathcal P}$ is a finite set of prime ideals of $R$ so that
no prime ideal in ${\mathcal P}$ contains the product $I_1I_2 \ldots I_r.$ 
Then there exists  a Rees-superficial element $a \in I_1$ for the ideals
$I_1, I_2, \ldots, I_r$ so that $a$ is not in any of the prime ideals in 
${\mathcal P}.$
\end{lemma}   
   
\begin{remark} 
It is clear that a nonzerodivisor in $I_1 \setminus I_1^2$  that is 
Rees-superficial for a set of ideals $(I_1, I_2, \ldots I_r)$ 
is also superficial.   Moreover  in  a Cohen-Macaulay local ring 
with infinite residue field, maximal Rees-superficial sequences that are also 
regular sequences exist for a set of $\m$-primary ideals, by  
Rees's basic lemma.  
\end{remark} 

For a function $f:\ZZ \lra \NN,$  put $\triangle f(n) = f(n)-f(n-1).$
\begin{pro} \label{diff} 

Let $(R,\m)$ be a local ring and $I$ an $\m$-primary ideal.  Let $a$
be a nonzerodivisor in $R$ which is Rees-superficial for $I$ and $\m.$  Let
``-'' denote residue classes in $\overline{R}=R/aR.$ Then for large
$n,$
$$
HF(F(\overline{I}),n)= \triangle HF(F(I),n).
$$
\end{pro} 

\begin{proof}
We have the exact sequence
$$
\CD 
0 @ > >> K_n  @ > >> I^n / \m I^n 
@ > \mu_a >> I^{n+1}/ \m I^{n+1} @ > >> C_n @ > >> 0, 
\endCD
$$
where $\mu_a( x + \m I^n) = ax + \m I^{n+1}$, $K_n=(\m I^{n+1} \;:\; a)\cap
I^n / \m I^n$ and $C_n= I^{n+1}/(aI^n+\m I^{n+1}).$ Since $a$ is
Rees-superficial for $\m$ and $I,$ $K_n=0$ for all large $n.$ Hence for
large $n,$ $ \triangle HF(F(I),n+1)= \mu(I^{n+1})- \mu(I^n)=
\ell(C_{n}).$ For all large $n$, 
\begin{eqnarray*}
HF(F(\bar{I}), n+1)
&=&\ell\left(I^{n+1}+aR/(\m I^{n+1}+aR)\right)\\
&=&\ell \left(I^{n+1}/(\m I^{n+1}+aR\cap I^{n+1})\right)\\
&=& \ell \left(I^{n+1}/(\m I^{n+1}+a I^n)\right)\\
&=& \mu(I^{n+1})-\mu(I^n)\\
&=& \triangle HF(F(I),n+1).
\end{eqnarray*}
\end{proof}
 
\begin{thm}
  Let $(R,\m)$ be a Cohen-Macaulay local ring of positive dimension $d.$  
  Let $I$ be an $\m$-primary ideal. 

\begin{enumerate}
  \item Let $a_1, a_2, \ldots, a_{d-1}\in I , x \in \m$ be a
    regular sequence in $R$ which is Rees-superficial sequence for 
    the multiset $(I^{[d-1]} |\m^{[1]}).$ 
    Then
    $$
    f_0(I)= e_{d-1}(\m|I) - \lim_{n \ra \infty} \ell \left( \frac{\m
    I^n } { xI^n+(a_1, a_2, \ldots, a_{d-1})\m I^{n-1} }\right).
    $$
  
  \item If $a_1, a_2, \ldots, a_d \in I$ is a regular sequence in $R$ which 
        is a Rees superficial sequence for
        for the multiset $(I^{[d-1]} |\m^{[1]}).$  Then  
    $$
    f_0(I)=e(I)- \lim_{n \ra \infty} \ell \left( \frac{\m I^n } {
    a_dI^n+(a_1, a_2, \ldots, a_{d-1})\m I^{n-1} }\right).
    $$
\end{enumerate}
\end{thm}
\begin{proof} 
(1) We induct on $d.$ Let $d=1.$  We need to prove that
$$
f_0(I)=e(\m)- \lim_{n \ra \infty} \ell \left( \frac{\m I^n } {
xI^n}\right).
$$
For all $n \in \NN,$
$
\ell(R/xR)+\ell(xR/xI^n)=\ell(R/\m I^n)+\ell(\m I^n/xI^n).
$
Since $x \in \m$ is superficial for $\m,$ $xR $ is a minimal reduction
of $\m.$ Therefore we get $\mu(I^n) = e(\m) - \ell(\m I^n/xI^n)$.
Hence by taking limits we get the desired formula.

Now suppose $d=2.$ Let $(a,x)$ be regular sequence which is 
a Rees-superficial sequence for $I,\m.$ Then 
$(a,x)$ is a joint reduction of the set $(I,\m),$ by Threorem \ref{jr}.  
And 
$e_1(\m |I)=\ell(R/(a,x))$ by \cite[Theorem 2.4(ii)]{re}. 
By the proof of Lemma 4.2 of
\cite{dv}, we have for all $n \geq 1,$ $$\triangle HF(F(I),n) =e_1(\m
|I)-\ell\left(\frac{\m I^n}{xI^n+a\m I^{n-1}}\right)+ \ell
\left(\frac{(\m I^{n-1}:x)\cap (I^n:a)} {I^{n-1}}\right).$$ 
Since $(a,x)$ is superficial for $(I,\m)$, $a$ is superficial for $I$
and is regular in $R$, we have, for large n, $I^n : a =
I^{n-1}$.
Since $HP(F(I),n)$ is a degree one polynomial, $\triangle
HF(F(I),n)=f_0(I)$ for large $n.$ This establishes the formula for
$d=2.$

Now suppose $d \geq 3.$ Put $\bar{R}=R/(a_1)$ and   $L=(a_2, a_3,
\ldots, a_{d-1}).$ By induction hypothesis and Proposition \ref{diff},
\begin{eqnarray*}
f_0(I)&=& f_0(\bar{I})\\
&=& e_{d-2}(\bar{\m}|\bar{I})- \lim_{n \ra \infty}
\ell \left(\frac{\m I^n+a_1R}{xI^n+ L \m I^{n-1}+a_1R}\right) \\
&= & e_{d-1}(\m|I)-\lim_{n \ra \infty}  \ell\left
( \frac{\m I^n+a_1R}{xI^n+L\m I^{n-1}+a_1R}\right)\\
&= & e_{d-1}(\m|I)- \lim_{n \ra \infty} \ell\left(\frac{\m I^n}{xI^n+ 
L \m I^{n-1}+\m I^n \cap a_1R } \right)\\
& = & e_{d-1}(\m|I)- \lim_{n \ra \infty}  
\ell\left( \frac{\m I^n}{xI^n+(a_1, \ldots, a_{d-1}) \m I^{n-1}}
\right).
\end{eqnarray*}
In the above equations we have used the fact that if  $a_1$ is
superficial for $\m$ and $I,$ then $e_{d-2}(\bar{\m}|\bar{I}) =
e_{d-1}(\m|I)$ by \cite[p.\ 118, line 3]{kv}.  This establishes the formula.

\noindent
(2) Replace $x$ by $a_d$ in the above argument.  
\end{proof}

We now obtain a sufficient condition for $f_0(I)=e_{d-1}(\m|I)$.

\begin{thm} \label{multrees}
Let $(R,\m)$ be a Cohen-Macaulay local ring of dimension $d$ and $I$
an $\m$-primary ideal. If $BF(r,s)=\ell(R/\m^rI^s) = BP(r,s)$ for all
$r, s \geq 0$, then 
$$
HS(F(I), t) = \frac{\sum_{j=0}^{d-1}(1-t)^{d-j-1}g(j)}{(1-t)^d}, 
$$
where $g(j) = \sum_{i=1}^{d-j}ie_{(i,j)}.$ In particular
$f_0(I)=e_{d-1}(\m|I).$ 
\end{thm}

\begin{proof}
For convenience, write $e_{(i,j)}
= e(i,j)$. Then
\begin{eqnarray*}
\mu(I^s) & = & \ell(R/\m I^s) - \ell(R/I^s) \\
& = & \sum_{i+j \leq d} e(i,j){1+i \choose i}{s+j \choose j} -
\sum_{i+j \leq d} e(i,j){i \choose i}{s+j \choose j} \\
& = & \sum_{i+j \leq d} ie(i,j) {s+j \choose j} \\
& = & \sum_{j=0}^{d-1}\left[\sum_{i=1}^{d-j} ie(i,j)\right]
{s+j \choose j}\\
& = & \sum_{j=0}^{d-1} g(j){s+j \choose j}.
\end{eqnarray*}
Hence we have
\begin{eqnarray*}
 HS(F(I),t) & = & \sum_{s\geq0}\mu(I^s)t^s =
 \sum_{s\geq0}\left[\sum_{j=0}^{d-1} g(j) {s+j \choose
 j}\right]t^s \\
 & = & \sum_{j=0}^{d-1}g(j)\left[\sum_{s\geq0}{s+j\choose
 j}t^s\right] \\
 & = & \sum_{j=0}^{d-1}\frac{g(j)}{(1-t)^{j+1}} \\
 & = & \frac{\sum_{j=0}^{d-1}(1-t)^{d-j-1}g(j)}{(1-t)^d}.
\end{eqnarray*}
Now put $t=1$ in the numerator of $HS(F(I),t)$ to get 
$f_0(I) = e_{d-1}(\m | I)$.  
\end{proof}

The next result  was communicated to us by E. Hyry.
\begin{cor} 
Let $(R,\m)$ be a Cohen-Macaulay local ring.  Let  the multi-Rees
algebra ${\mathcal R} := R[\m t_1, It_2]$ be Cohen-Macaulay. Then
$f_0(I)=e_{d-1}(\m | I).$ If $d = 2,$ then $R$ and $F(I)$ are
Cohen-Macaulay with minimal multiplicity.
\end{cor} 

\begin{proof}
If  $\mathcal R$ is Cohen-Macaulay, then by
the proof of \cite[Theorem 6.1]{h}, $\ell(R/\m^rI^s) = BP(r,s)$ 
for all $r, s \geq 0$.
Therefore $f_0(I)=e_{d-1}(\m | I),$ by Theorem \ref{multrees}. When
$d=2$ and  $\mathcal R$ is Cohen-Macaulay, then, by \cite[Corollary 3.5]{h2}, 
$R[\m t_1]$ and $R[It_2]$ are Cohen-Macaulay and hence $r(I) \leq 1$ and
$r(\m) \leq 1$ by \cite[Remark 3.10]{gs}. 
Thus $R$ is Cohen-Macaulay with minimal
multiplicity. Since $r(I) \leq 1$, by \cite[Theorem 1 and 7]{shah} 
$F(I)$ is Cohen-Macaulay and $f_0(I) = \mu(I) -1$. Hence  $F(I)$ 
has minimal multiplicity.
\end{proof}

\section{\bf Cohen-Macaulay fiber cones}
In this section we use the multiplicity formula for fiber cones to
detect their Cohen-Macaulay property. We begin by recovering Corollary
2.5 of \cite{drv} in a simpler way.

\begin{pro}
Let $(R,\m)$ be a $d$-dimensional Cohen-Macaulay local ring and $I$ an
$\m$-primary ideal of minimal mixed multiplicity.  Then $F(I)$ is
Cohen-Macaulay if and only if $r(I) \leq 1.$
\end{pro}

\begin{proof}
Let $J$ be any minimal reduction of $I.$ Then $F(I)$ is Cohen-Macaulay
if and only if $f_0(I)=\ell (F(I)/JF(I)).$ Since
$$
\frac{F(I)}{JF(I)}=
               \frac{R}{\m} \oplus \frac{I}{J+\m I} \oplus 
              \left( \bigoplus_{n=2}^{\infty}\frac{I^n}{JI^{n-1}+\m I^n}
              \right),  
$$
and $\ell\left(I/(J+\m I)\right) = \mu(I) -d,$ we have 
$$\ell\left(F(I)/JF(I))\right)= 1+ \mu(I) -d + 
\sum_{n=2}^{\infty}\ell\left(\frac{I^n}{\m I^n+JI^{n-1}}\right).$$
Thus $F(I)$ is Cohen-Macaulay if and only if for any Rees-superficial
sequence $a_1, a_2, \ldots, a_{d-1}, x,$  where  
$a_1, a_2, \ldots, a_{d-1} \in I , x \in \m$ 
\begin{eqnarray*}
f_0(I)&=& e_{d-1}(\m|I)- \lim_{n \ra \infty}\ell 
\left(\frac{\m I^n}{xI^n+(a_1, a_2, \ldots, a_{d-1})\m I^{n-1}}\right) \\
&=& \mu(I)-d +1 +
\sum_{n=2}^{\infty}\ell\left(\frac{I^n}{\m I^n+JI^{n-1}}\right).
\end{eqnarray*}
By the proof of Proposition 2.4 of \cite{drv}, $I$ has minimal mixed
multiplicity if and only if  $I^n\m=xI^n+ (a_1, a_2, \ldots,
a_{d-1})\m I^{n-1}$ for all $n \geq 1.$ Thus $F(I)$ is
Cohen-Macaulay if and only if $I^2=JI.$
\end{proof}
In the next result we improve the Corollary 1.4 of \cite{dv} by
removing the hypothesis of almost maximal depth for the associated
graded ring of $I.$

\begin{pro} \label{ammm}
Let $(R,\m)$ be a Cohen-Macaulay local ring with infinite residue
field. Let $I$ be an $\m$-primary ideal with almost minimal mixed
multiplicity. Then 
\begin{enumerate} 
\item Either $f_0(I)=e_{d-1}(\m | I)$ or
  $f_0(I)=e_{d-1}(\m | I)-1.$

\item Let $f_0(I) =e_{d-1}(\m|I).$ Then $F(I)$ is Cohen-Macaulay if
  and only if $r(I)=2$ and $\ell \left(I^2/(JI+\m I^2)  \right)=1.$

\item  Let $f_0(I) =e_{d-1}(\m|I)-1.$ Then $F(I)$ is Cohen-Macaulay if
  and only if $r(I) \leq 1.$
\end{enumerate}
\end{pro}

\begin{proof} 
(1) Let $a_1, a_2, \ldots, a_{d-1} \in  I  ,x \in \m$ be a Rees-superficial
sequence for $I$ and $\m.$ Put $L=(a_1, a_2, \ldots, a_{d-1})$ and
$\alpha=\lim_{n \ra \infty}\ell\left(\m I^n/(xI^n+L\m I^{n-1})
\right).$ Since $I$ has almost minimal mixed multiplicity, 
$\ell \left(\m I^n/(xI^n+L \m I^{n-1})\right) \leq 1$ for
all $n$ by \cite[Lemma 2.2]{dv}. Hence $\alpha=0$ or  $1.$ This proves (1).
\vskip 2mm
\noindent
(2) By the computations in the above result,
  $F(I)$ is Cohen-Macaulay 
\begin{align*} 
       \iff & e_{d-1}(\m|I)-\alpha = \mu(I)-d+1 +
\sum_{n=2}^{\infty}\ell\left(\frac{I^n}{\m I^n+JI^{n-1}}\right),
 \\
 \iff &   1-\alpha= \sum_{n=2}^{\infty}\ell
        \left(\frac{I^n}{\m I^n+JI^{n-1}}\right).
\end{align*} 

Let $f_0(I) =e_{d-1}(\m|I).$ Then $\alpha=0.$ Thus
 
\begin{align*}
F(I) \ \text{is Cohen-Macaulay} \ &
\iff \sum_{n=2}^{\infty}\ell\left(I^n/(\m I^n+JI^{n-1})\right)=1, \\
 &\iff \ r(I)= 2 \ \text{and} \ \ell(I^2/(JI+\m I^2))=1.
\end{align*}

(3) Let
$f_0(I) =e_{d-1}(\m|I)-1.$ Hence $\alpha=1.$  Thus $F(I)$ is
Cohen-Macaulay if and only if 
$$\sum_{n=2}^{\infty}\ell\left(I^n/\m I^n+JI^{n-1}\right)=0.$$ This holds
if and only if $I^2=JI.$
\end{proof}

G. Valla raised, in a personal communication, a question regarding
the Cohen-Macaulay property of fiber cones of Sally ideals. In 
 Example 6.1  we show that $F(I)$ need not be Cohen-Macaulay, even
if $G(I)$ is Cohen-Macaulay. First we characterize Cohen-Macaulay
fiber cones of Sally ideals in dimension one.

\begin{thm} \label{dim1}
Let $(R,\m)$ be a $1$-dimensional Cohen-Macaulay local ring, $I$ a
Sally  ideal and $ J = (x)$ a minimal reduction of $I$, with reduction
number $r$. Then the following are equivalent:
  \begin{enumerate}
    \item[(i)] $F(I)$ is Cohen-Macaulay.
    \item [(ii)] $HS(F(I),t)= (1+(\mu(I)-1)t+t^2+t^3 +\cdots+t^r)/(1-t).$
\item[(iii)] $\mu(I^k) = \mu(I) + k - 1,$ 
                 \mbox{ for } $2 \leq k \leq r$.    
    \item[(iv)] $\mu(I^2) = \mu(I) + 1$.
\item[(v)] $\m I^2 = \m JI$.

  \end{enumerate}
\end{thm}

\begin{proof}
(i)$ \Rightarrow$ (ii): Let $F(I)$ be Cohen-Macaulay. Then by \cite[Theorem 2.1]{drv},
$$HS(F(I),t)=\frac{1+(\mu(I)-1)t+
 \sum_{i=2}^r\ell(I^i/(JI^{i-1}+\m I^i)t^i}{(1-t)}
.$$
Since $I$ is a Sally ideal $\m I^n \subset JI^{n-1}$ for all $n \geq 2$
and $\ell(I^n/JI^{n-1})=1$ for all $n=2, 3, \ldots, r.$ Hence (ii) follows.

\noindent
(ii)$\Rightarrow $(iii): From the formula for the Hilbert series of $F(I)$, 
we obtain the equation $\mu(I^k)=\mu(I)+k-1$ for $k=1, 2, \ldots, r.$

\noindent
(iii)$\Rightarrow $(iv): Put $k=2.$

\noindent
(iv)$\Rightarrow $(v): For $n \geq 1$ 
we have  following exact sequence:
$$
0 \longrightarrow \frac{(\m I^{n+1} : x) \cap I^n}{\m I^n}
\longrightarrow \frac{I^n}{\m I^n}
{\overset{\phi_x}{\longrightarrow}} \frac{I^{n+1}}{\m I^{n+1}}
\longrightarrow \frac{I^{n+1}}{xI^n} \longrightarrow 0,
\eqno(1)
$$
and the isomorphism:
$$
\frac{(\m I^{n+1} : x) \cap I^n}{\m I^n} \cong \frac{xI^n \cap
\m I^{n+1}}{x\m I^n}. \eqno(2)
$$
Assume $\mu(I^2)
= \mu(I) + 1$. Then, from the exact sequence (1) and the isomorphism
(2) for $n = 1$, we get $\m I^2 = x\m I$.
\vskip 2mm
\noindent

\noindent
(v) $\Rightarrow$ (i):  Consider the function
$H_{\m}(I,n) := \ell(R/\m I^n)$ and write the corresponding polynomial as: 
$$
P_{\m}(I,n) = \sum_{i=0}^d (-1)^ig_i(I){n+d-i-1 \choose d}.
$$
Then by \cite[Theorem 5.3]{jv}, 
$$g_1(I) =  \sum_{n\geq1} \ell(\m I^n/x\m I^{n-1}) - 1. $$
Since $\m I^2 = x\m I$ we get
$$ g_1(I) = \ell(\m I /x\m)-1. \eqno(3)$$

We know  $F(I)$ is
Cohen-Macaulay if and only if $g_1(I) = \sum_{n\geq1}\ell(\m I^n +
xI^{n-1}/xI^{n-1}) - 1$ by  \cite[Theorem 4.3]{jv}. 
Since $\ell(I^2/xI) = 1,$ $\m I^n \subset
xI^{n-1}$ for all $n \geq 2$.
Therefore, by (3),
$$
 \sum_{n\geq1}\ell(\m I^n + xI^{n-1}/xI^{n-1}) -1= 
 \ell(\m I + xA/xA) -1 = \ell(\m I/x\m) -1 
  =  g_1(I).
$$ 

Hence  $F(I)$ is Cohen-Macaulay. 
\end{proof}

Now we characterize the Cohen-Macaulayness of $F(I)$ in higher
dimensions.

\begin{thm}
Let $(R,\m)$ be a Cohen-Macaulay local ring of dimension $d \geq 1$, 
$I$ a Sally ideal with a minimal reduction $J$. Then the following are
equivalent:

\begin{enumerate}
         \item[(i)] $F(I)$ is Cohen-Macaulay.
         \item[(ii)] $\m I^2 = \m JI$.
         \item[(iii)] The Hilbert series of $F(I)$ is given by  $$HS(F(I),t) =
                  \frac{1+(\mu(I)-d)t+t^2+\cdots+t^r}{(1-t)^d}.$$
         \item[(iv)] $f_0(I) = \mu(I) - d + r$.
\end{enumerate}
\end{thm}

\begin{proof} We apply induction on $d.$ 
We have proved the theorem for $d=1.$  Now let $d \geq 2.$ \\
(i) $\Rightarrow$ (ii): Since $I$ is a Sally ideal, by \cite[Corollary 1.7]{r}, 
depth $G(I) \geq d-1.$ Hence we can choose an $x \in J$ such that 
$x^*$ is regular in
$G(I)$ and $x^o$ is regular in $F(I)$. Then $F(I/(x))\cong F(I)/(x^o)$ 
is Cohen-Macaulay. By induction, $\bar{\m}\bar{I}^2 =
\bar{\m}\bar{J}\bar{I}$. Therefore $\m I^2 = \m JI + (x) \cap \m I^2$.
Since $x^*$ is regular in $G(I)$ and $x^o$ regular in $F(I)$, $(x)
\cap \m I^2 = x\m I$, by \cite[Theorem 1.1]{cz}. Therefore $\m I^2 = \m JI$.
\vskip 2mm
\noindent
(ii) $\Rightarrow$ (i): For $x \in J$, such
that $x^*$ is regular in $G(I)$ and $x^o$ superficial for $F(I)$,
let ``--'' denote ``modulo $(x)$''. Then $\bar{\m}\bar{I}^2 =
\bar{\m}\bar{J}\bar{I}$. By induction, $F(\bar{I})$ is Cohen-Macaulay.
By ``Sally machine'' \cite[Lemma 2.7]{jv}, $x^o$ is regular in 
$F(I)$ and hence $F(I)$ is
Cohen-Macaulay.
\vskip 2mm
\noindent (i) $\Rightarrow$ (iii): 
 Since $F(I)$ is Cohen-Macaulay,
 $$HS(F(I),t) = \frac{HS(F(I)/JF(I),t)}{(1-t)^d}.$$ 
Since $\m I^2 \subset JI$, we have
$\ell(I^n/\m I^n+JI^{n-1}) = \ell(I^n/JI^{n-1}) = 1$ for all $n
= 2, \ldots, r$. Therefore the Hilbert series of $F(I)$ is
$$
HS(F(I),t) = \frac{1+(\mu(I)-d)t+t^2+\cdots+t^r}{(1-t)^d}. 
$$
\vskip 2mm
\noindent(iii) $\Rightarrow$ (iv): 
The assertion follows directly from the fact that if $HS(F(I),t) =
h(t)/(1-t)^d$, then $f_0(I) = h(1)$.
\vskip 2mm
\noindent (iv) $\Rightarrow$ (i): 
Since $f_0(I) = \mu(I)-d+r$, we have
\begin{eqnarray*}
  1 + \sum_{n=1}^r \ell\left(\frac{I^n}{JI^{n-1}+\m
  I^n}\right) & = & 1 + \ell(I/\m I+J) + 
  \sum_{n=2}^r\ell(I^n/JI^{n-1})\\
  & = & 1 + \mu(I) - \ell(\m I+J/\m I) + r - 1 \\
  & = & \mu(I) - \ell(J/\m J) + r \\
  & = & \mu(I) - d + r \\
  & = & f_0(I).
\end{eqnarray*}
Therefore, by \cite[Theorem 2.1]{drv}, $F(I)$ is Cohen-Macaulay.
\end{proof}
 
\section{\bf Gorenstein  fiber cones }

{\em Throughout this section and the next we will assume, unless
otherwise stated, $(R,\m)$ is a Cohen-Macaulay  local ring of
dimension $d$ with infinite residue field, $I$ is an $\m$-primary
ideal and $F(I)$ is Cohen-Macaulay.}

 In this section we study Gorenstein property of $F(I)$ for several 
classes of ideals. We do this by keeping reduction numbers in mind. 
 It is clear that if $r(I) = 0,$ then $F(I)$ is a  polynomial ring. 
Thus we may begin with  the case when $r(I) =1$. 
In this case $F(I)$ is Cohen-Macaulay   \cite[Theorem 1]{shah}.
\begin{pro}
\label{fgore} Assume  $r(I) =1$.  If $F(I)$ is
Gorenstein, then $\mu(I) = d + 1$. 
\end{pro}
\begin{proof}
The Hilbert series of $F(I)$ is $(1 + (\mu(I) -d)t)/(1-t)^d$ by \cite[Theorem 2.1]{drv}.
By a theorem of Macaulay cf. \cite[Theorem 4.1]{stan}, the $h$-vector of a standard Gorenstein  graded 
$k$-algebra, where $k$ is a field, is symmetric. Hence  $\mu(I) - d =1$. 
\end{proof}

\begin{remark}
The symmetry of the $h$-vector does not imply that the Fiber cone is Gorenstein even when it is
Cohen-Macaulay, see Example 6.2
\end{remark}

In general we have the following:
\begin{pro}
\label{mugore} If $\mu(I) = d+1,$ then $F(I)$ is a hypersurface.  
\end{pro}
\begin{proof}
 Set $I = (u_1,\dots,u_{d+1})$.
Consider the map $
\phi \colon k[X_1,\ldots,X_{d+1}] \longrightarrow F(I)$ given by
             $\phi(X_i)= u_i + \m I $ for $ i =1,\ldots,d+1.$
Clearly $\phi$ is surjective and  $\ker (\phi)$ is 
a height one ideal of $S = k[X_1,\ldots,X_{d+1}]$. 
Since $ S/\ker(\phi) \cong F(I)$ is Cohen-Macaulay, $\ker (\phi)$ is
a height one unmixed ideal. Since $S$ is a UFD, $I$ is principal. 
Therefore $F(I)$ is a hypersurface ring.
\end{proof}

\begin{remark}
One of the surprising results in our investigations has been the
following: If $F(I)$ is Gorenstein, then $\mu(I)$ is forced. When
$r(I) = 1$ this is done in Proposition \ref{fgore}. When $r(I) = 2$
and $I$ has almost minimal multiplicity we get an upper bound on
$\mu(I)$, see Corollary \ref{gorbound}.
\end{remark}

 Let $J  = (x_1,\ldots,x_d)$ be a minimal reduction of $I$. Since $F(I)$ 
 is Cohen-Macaulay, it easily follows that   the  reduction number of $I$ with 
respect to  $J$ is the degree of the $h$-polynomial of $F(I)$. We will use 
this fact implicitly in all subsequent discussions.
\begin{pro}
Set $r = r(I) $, the reduction number of $I$ and let  
$J$ be a minimal reduction of $I$. Then
$$
\socle F(I)/JF(I)  
\cong  \bigoplus_{n = 1}^{r-1}
 \frac{ \left(I^{n+1}\m + JI^n \colon I\right)
 \cap I^n}{\left(I^{n}\m + JI^{n-1}\right) }
 \bigoplus \frac{I^r}{\m I^r + JI^{r-1}}.
$$
\end{pro}
\begin{proof}
Let $J = (x_1,\ldots,x_d)$ be a minimal reduction of $I$. 
Since  $x_{1}^{\circ},\ldots, x_{d}^{\circ}$ is a regular sequence for 
$F(I),$ we have  
$$
\socle F(I)/JF(I) \cong \socle F(I)/(x_{1}^{\circ},\ldots, 
x_{d}^{\circ})F(I).
$$
Since  $S : = F(I)/JF(I)$   is a standard graded $k$-algebra with $k=
S_0,$ a field, we have $\socle S = (0 \colon_S \  S_+) = (0 \colon_S \
S_1)$, where $S_+ = \oplus_{n\geq1}S_n$.  Notice that $S_1 = I/(\m I +
J)$. An easy computation yields the result.
\end{proof}

If $r(I) = 2$, then we have the following ideal-theoretic condition to
check the Gorenstein property of $F(I)$.
\begin{cor}
\label{red2gore}

Let $ r(I) = 2 $ and let $J$ be a minimal reduction of $I$.
Then $F(I)$ is Gorenstein if and only if
$$
(I^2\m + JI \colon I)
 \cap I = \m I + J \quad \text{and} \quad 
 \ell\left(\frac{I^2}{\m I^2 + JI} \right) = 1.  
$$
\end{cor}
By using  Proposition \ref{ammm}, we get that, 
the  fiber cone $F(I)$ of  an ideal
$I,$  with $r(I)=2$ having almost minimal mixed multiplicity and
$f_0(I)=e_{d-1}(\m |I),$ is Gorenstein if and only if $(I^2\m + JI
\colon I) \cap I = \m I + J .$  If $r(I) \geq 3$ and $F(I)$ is
Gorenstein, then the symmetry of $h$-vector  yields the following
\begin{pro}
\label{red3gore}
Let $r =  r(I)\geq 3 $ and $J$ be a minimal reduction of $I$.
If $F(I)$ is Gorenstein, then 
$$
\mu(I) = d + \ell\left(\frac{ I^{r-1} }{\m I^{r-1} + J I^{r-2}} \right). 
$$
\end{pro}
\begin{proof}
Note that  $S : = F(I)/JF(I)$ is  a standard graded Gorenstein ring of 
dimension zero and 
$$
S_{r-1} = \frac{I^{r-1}}{\m I^{r-1} + JI^{r-2}} \quad 
\text{and} \quad S_{1} =  \frac{I}{\m I +J}.
$$
If $F(I)$ is Gorenstein, then the $h$-vector  of $F(I)$ is symmetric. 
Since $r \geq 3$ we have
$\ell (S_1) = \ell(S_{r-1})$. Since $\ell(S_1) = \mu(I) -d$ we get the result.
\end{proof}
As an easy consequence  we have:
\begin{cor}
\label{sallygore}
Let $I$ be a Sally ideal with $ r(I) \geq 3 $. If $F(I)$ is
Gorenstein, then
$\mu(I) = d+1$
\end{cor}
\begin{proof}
Let $J$ be a reduction of $I$. Then $I^2 \neq JI$. We have 
$\ell(I^{n+1}/JI^n) \leq 1$ for all $n \geq 1$. 
Notice that $\m I^{n+1} \subseteq JI^{n}$ for all $n \geq 1$.
Therefore we have
$$
 \frac{I^{r-1}}{\m I^{r-1} + JI^{r-2}}  =     \frac{I^{r-1}}{ JI^{r-2}}.
$$
Also $\ell( I^{r-1}/ JI^{r-2}) = 1$. The result follows from Proposition 
\ref{red3gore}.
\end{proof}
The result above does not hold if $r(I) \leq 2$. 
Consider the following example, discussed in \cite{s3}.  Let $e > 3$
be a positive integer. Set $R = k[\![t^e, t^{e+1}, \ldots,
t^{2e-2}]\!]$, where $k$ is a field. Since the numerical semigroup
generated by $\{e, e+1, \ldots, 2e-2\}$ is symmetric with conductor
$2e$, $R$ is Gorenstein. Let $\m$ denote the maximal ideal of $R$.
Then $\mu(\m) = e(R)+d-2 = e - 1 > d+1 = 2,$ where $e(R)$ denotes the
multiplicity of $R$.  By the proof of Theorem
3.4 of \cite{s3}, $\m^3 = J\m^2$ for any minimal reduction $J$ of $\m$
and it follows from Theorem 3.4 that $G(\m) = F(\m)$ is Gorenstein.

\section{\bf Gorenstein Fiber cones of ideals of almost minimal  multiplicity}
In this section we consider the Gorenstein property of fiber cones of
ideals of almost minimal multiplicity. Recall that an $\m$-primary
ideal $I$ in a Cohen-Macaulay local ring $(R,\m)$ is said to have
minimal multiplicity (resp. almost minimal multiplicity) if for any
minimal reduction $J$ of $I,$ $\m I = \m J$ (resp.  $\ell(I\m/J\m) =
1$). Such ideals have been studied in \cite{g1}, \cite{g2} and
\cite{jv2}.  

In addition to the hypotheses stated in the beginning of the previous
section we further assume $G(I)$ is Cohen-Macaulay. Since $I^3
\subseteq J$, we get by the  Valabrega-Valla criterion that  $I^3 =
JI^2$. So $r(I) \leq 2$.  Since we  have already considered the  case
$r(I) = 1$, we assume $r(I) = 2$.

Let $J$ be a minimal reduction of $I$. Set $J = (x_1,\ldots,x_d)$.  If
$G(I)$ is Cohen-Macaulay, then $x_{1}^{*},\ldots,x_{d}^{*}$ is a
$G(I)$-regular sequence. Since $ F(I)$ is Cohen-Macaulay, we also have
that $x_{1}^{\circ},\ldots,x_{d}^{\circ}$ is an  $F(I)$-regular
sequence.

\begin{notation}\label{notation}
Set $(B,\n) = (A/J, \m / J)$, $K = I/J$.
We have
\[
\frac{F(I)}{(x_{1}^{\circ},\ldots,x_{d}^{\circ})F(I)} \cong F(K).
\]
It follows that  $F(I)$ is Gorenstein if and only if  \ $F(K)$ is Gorenstein.
\end{notation}

\begin{blank}
\label{notar}
 Notice 
\begin{enumerate}
\item
$\n K \cong k$.
\item
$\n^2K = 0$ and so $K^3 = 0$ and $\n K^2 = 0$.
\item
$0 \neq K^2 \subseteq \n K$. So $K^2 = \n K$.
\end{enumerate}
\end{blank}
\begin{remark}
If $I$ has almost minimal multiplicity with $r(I) \geq 2$ and $I^2
\cap J = JI$, then $I$ is a Sally ideal. To see this, note that $K^2
\cong I^2/J\cap I^2 = I^2/JI$ and from \ref{notar}(1) and (3), it
follows that $\ell(I^2/JI) = 1$. In particular, if $G(I)$ is
Cohen-Macaulay, then $I$ is a Sally ideal.
\end{remark}

\begin{remark}
Since $r(I) = 2$, symmetry of the $h$-vector of Hilbert series of $F(I)$ 
does not  help us in estimating $\mu(I)$. To find conditions on  $\mu(I)$ 
we need the following different  criterion.
\end{remark}
\begin{pro}
\label{ammt}
 Assume that $I$ has almost minimal multiplicity. 
Set $W = I\cap (\m J \colon I)$. 
Then $F(I)$ is Gorenstein  if and only if $W = \m I + J$. 
\end{pro}
\begin{proof}
Since $r(I) = 2$ we can use Corollary \ref{red2gore}.  We keep
notation as in \ref{notation}. Note that $K^2 /\n K^2 = 
(I^2 + J)/(\m I^2 +J)$. Then we have
$$
 \frac{I^2 + J}{\m I^2 +J} \cong \frac{I^2}{(\m
I^2 + J) \cap I^2} = \frac{I^2}{\m I^2 + JI}.
$$
Thus $\ell (I^2 /(\m I^2 + JI)) = 1$.
Set $E = (I^2\m + JI \colon I) \cap I$.  Using Corollary \ref{red2gore} 
we have that $F(I)$ is Gorenstein if and only if $E = \m I + J$.

We now prove that $E = W$.  Since $\ell(\m I/\m J) = 1$, we have that
$\m I^2 \subseteq J\m$.  It follows that $W \supseteq E$. Conversely,
let $t \in W$. Then $tI \in \m J$. In particular $tI \subseteq J \cap
I^2 = JI$. So $t \in E$.
Therefore $E = W$.
\end{proof}

\begin{remark}
The hypothesis $G(I)$ Cohen-Macaulay is essential in the Proposition above.
See Example 6.4.
\end{remark}

\begin{blank}
{\rm  We define a  bilinear form,
\begin{align*}
\phi \colon \frac{K}{\n K} \times \frac{\n}{\n^2K \colon K} &
\longrightarrow \frac{\n K}{\n^2 K} \\ 
\left( a + \n K , b + \n^2K\colon K \right) &\mapsto  ab + \n^2 K.
\end{align*}
Notice that with our hypothesis $\n^2 K = 0$. 
It is straightforward to check that $\phi $ is well defined. 
It is also clear that
$\phi $ is non-degenerate with respect to $\n/(\n^2K \colon K)$.
}
\end{blank}

\begin{lemma} \label{ndgen}
 Assume that $I$ has almost minimal multiplicity.
If $W = \m I + J,$ then $\phi$ is non-degenerate with respect to  $K/\n K$.
\end{lemma}
\begin{proof}
Suppose there exists $a \in I$ (so $\ol{a} \in K$) such that
$$
 \phi \left(\ol{ a} + \n K , b + \n^2K \colon K \right) 
=  ab + \n^2 K = \ol{0} \quad \text{for each $b + \n^2 \in  \n/\n^2K:K.$}
$$
Thus $\ol{a} \n \subseteq \n^2 K = \{ \ol{0} \}$. So $ a \m \subseteq J$.
In particular $aI \subseteq J$. So $aI \subseteq J \cap I^2 \subseteq J\m$. 
 Thus $a \in \m J \colon I \cap I$. Thus $a \in W = \m I + J$. 
So $\ol{a}=0$. Therefore  $\phi$ is non-degenerate with respect to  $K/\n K$.
\end{proof}

\begin{cor}
\label{gorbound}
Assume that $I$ has almost minimal multiplicity and 
let $J$ be a minimal reduction of $I$.
If $F(I)$ be Gorenstein, then
$\ell(J:I/J) = \ell (R/I)$ and  $\mu(I) \leq  \mu(\m) + d$.
\end{cor}
\begin{proof}
Let $x_1,\ldots,x_d$ be a superficial sequence in $R$ with respect to
$I$. Set $J = (x_1,\ldots,x_d)$ and $(B,\n) = (R/J, \m/J)$ and $K =
I/J$.  Since $F(I)$ is Gorenstein, by Proposition \ref{ammt}  $W = \m
I + J$.  Therefore, as in the  Lemma \ref{ndgen}, $\phi$ is
non-degenerate both on the left and the right.  It
follows that $ \ell( K/\n K) =  \ell (\n/(0:K))$(notice that we used 
$\n^2K =  0$).
 This immediately
yields $$\mu(I) - d  = \mu(K)= \ell (\n/(0:K))  \leq \mu(\n) \leq  \mu(\m).$$
Notice $0:K \cong (J:I)/J$. We have 
$$\mu(I)-d  = \mu(K) = \ell (\n/(0:K)) =\ell (\m /(J:I)),$$ 
and
$$\ell (\m /(J:I)) =\ell(\m /J) - \ell(J:I/J) = e(I) - 1 - \ell(J:I/J).$$
So $\mu(I) = d+  e(I) - 1 - \ell(J:I/J) $. Since $I$ has 
almost minimal mixed multiplicity we have
$\mu(I) =  d+  e(I) - 1 - \ell(R/I)$. So we get 
$ \ell(J:I/J) = \ell (R/I)$.
\end{proof}

\begin{remark}{\em
The corollary above has an interesting connection with the Koszul
homologies, $H_{i}(I,R)$, of $I$. Let $\mu(I) = n$ and let $J = (x_1,
\ldots,x_d)$ be a minimal reduction of $I$. Since $R$ is
Cohen-Macaulay and $I$ is $\m$-primary we have (see \cite[1.6.16 and
1.6.17]{bh}) 
\begin{eqnarray}
& & \left\{
\begin{array}{ll}
H_i(I,R) = 0 \qquad \text{for}\  & i > n - d. \\
H_{n-d}(I,R) \cong \frac{J\colon I}{J}, \qquad H_{0}(I,R) &= \frac{R}{I}.
\end{array} \right.
\end{eqnarray}
Thus (with the hypotheses as in Corollary \ref{gorbound}), we get that
the zeroth Koszul homology and the last non-vanishing Koszul homology
of $I$ have the same length. Note that this property is true
if $R$ is Gorenstein and $I$ is any $\m$-primary ideal.
We have not been able to find an example where
$R$ is not Gorenstein, $G(I)$ is Cohen-Macaulay and $F(I)$ is Gorenstein
and $(1)$ is satisfied. However, we do not believe that
if for an $\m$-primary ideal $I$ in a Cohen-Macaulay local ring $R$ 
such that $G(I)$ is Cohen-Macaulay and $F(I)$ is Gorenstein, then
$R$ is Gorenstein.
}
\end{remark}
\section{\bf  Examples}
We end this article by presenting a few  examples to illustrate our
results. Let $k$ denote a field. The computations have been
performed in CoCoA, \cite{co}.

\begin{example}{\rm 
Let $A = k[\![t^6, t^{11}, t^{15}, t^{31}]\!]$, $I = (t^6, t^{11},
t^{31})$ and $J = (t^6)$. Then, it can easily be verified that 
$\ell(I^2/JI) = 1$ and $I^3 = JI^2$. Since $I^2 \cap J = JI$, 
$G(I)$ is Cohen-Macaulay. It can also be seen that $t^{37} \in \m
I^2$, but $t^{37} \notin \m JI$. Therefore $F(I)$ is
not Cohen-Macaulay. This example shows that the length condition
$\ell(I^2/JI)=1$ by itself need not force the fiber cone to be
Cohen-Macaulay, even if $ G(I) $ is Cohen-Macaulay.
}
\end{example}

\begin{example}{\rm
Here we give an example of an ideal $I$ such that $F(I)$ and
$G(I) $ are Cohen-Macaulay, the numerator of the Hilbert series is
symmetric, but $F(I)$ is not Gorenstein. Consider $A = k[\![t^7,
t^{15}, t^{17}, t^{33}]\!]$, $I = (t^7, t^{17}, t^{33})$ and $J =
(t^7)$. Then $\ell(I^2/JI) = 1, \; I^3 = JI^2, \; I^2 \cap J = JI$
and $\m I^2 = \m JI$. Therefore, in this case $ G(I) $ and $F(I)$
are Cohen-Macaulay. Hence, the Hilbert series of $F(I)$ is
$$
HS(F(I),t) = \frac{1+(\mu(I)-d)t+t^2}{(1-t)} = \frac{1+2t+t^2}{(1-t)}.
$$
It can also be seen that $t^{33} \in (\m I^2+JI : I) \cap I$, but
$t^{33} \notin J + \m I$.  Therefore, the numerator of the Hilbert
series is symmetric, but $F(I)$ is not Gorenstein.
}
\end{example}

\begin{example}{\rm 
Now we give an example of an ideal with Gorenstein fiber cone which is
not a hypersurface. Let $A = k[\![t^4, t^5, t^6, t^7]\!]$, $I = (t^4,
t^5, y^6)$ and $J = (t^4)$. Then, $\ell(I^2/JI) = 1, \; I^3 = JI^2$
and $\m I = \m J$. Therefore, $F(I)$ is Cohen-Macaulay. Since $t^{11}
\in I^2 \cap J$ and $t^{11} \notin JI, \; G(I) $ is not
Cohen-Macaulay. It can also be easily checked that $(\m I^2 + JI : I)
\cap I = J + \m I$. Hence $F(I)$ is Gorenstein.
}
\end{example}

\begin{example}{\rm
Let $A = k[\![x,y]\!], \; I = (x^3, x^2y, y^3)$ and $J = (x^3, y^3)$.
Then $\ell(I^2/JI) = 1, \; I^3 = JI, \; \ell(\m I/\m J) = 1$ and
$\m I^2 = \m JI$. Since $\mu(I) = d+1$, $F(I)$ is a hypersurface. It
can easily be seen that $x^4y^2 \in I^2 \cap J$, but $x^4y^2 \notin
JI$. Hence $ G(I) $ is not Cohen-Macaulay. It is easily checked that
$(\m J : I) \cap I \neq J +\m I$, even if $F(I)$ is Gorenstein. This
shows that the assumption on the Cohen-Macaulayness of $ G(I) $ in
Theorem \ref{ammt} is necessary.
}
\end{example}

\begin{example}{\rm 
Let $A = k[\![x,y,z]\!], \; I = (x^3, y^3, z^3, xy, yz, zx)$ and $J =
(x^3+yz, y^3+z^3+xz, xz+xy)$. It is been shown in \cite{drv} that $I$
has minimal mixed multiplicity. It can be seen that $I^2 = JI$. 
Therefore both $G(I) $ and $F(I)$ are Cohen-Macaulay. It can also be
seen that $\ell(\m I/\m J) = 1$. Since $z^3 \in \m J : I \cap I$ and
$z^3 \notin \m I + J$, $F(I)$ is not Gorenstein, by Proposition
\ref{ammt}. The Hilbert series of the fiber cone is $HS(F(I),t) =
1+3t/(1-t)^3$, which also shows that the fiber cone is not Gorenstein.
But, we have the equalities, $\ell(J:I/J) = \ell(R/I) = 7$ and $\mu(I)
= 6 = \mu(\m) + d$. This shows that the converse of Corollary 
\ref{gorbound} is not true.
}
\end{example}

\end{document}